%Paper: dg-ga/9411006
%From: huebschm@mpim-bonn.mpg.de (Johannes Huebschmann)
%Date: Tue, 22 Nov 1994 22:28:31 --100

%&amstex
%This file can be processed with AMSTeX 2.1
\documentstyle{amsppt}
\magnification=1200
\hoffset=-0.5pc
\nologo
\vsize=57.2truepc
\hsize=38.5truepc
\spaceskip=.5em plus.25em minus.20em
\define\abramars{1}
\define\armgojen{2}
\define\armamonc{3}
\define\atibottw{4}
 \define\ebinone{5}
 \define\golmiltw{6}
\define\singulat{7}
\define\singulth{8}
\define\topology{9}
 \define\smooth{10}
  \define\direct{11}
 \define\locpois{12}
 \define\modusym{13}
\define\huebjeff{14}
\define\marswein{15}
\define\mittvial{16}
\define\nararama{17}
\define\narashed{18}
\define\isingeth{19}
\define\sjamaone{20}
\define\sjamlerm{21}
\define\weinstwe{22}
\noindent
dg-ga/9411006
\newline\noindent
Math. Z., to appear.
\bigskip
\topmatter
\title The singularities of Yang-Mills connections for\\
       bundles on  a surface. I. The local model
\endtitle
\author Johannes Huebschmann{\dag}
\endauthor
\affil
Max Planck Institut f\"ur Mathematik
\\
Gottfried Claren-Str. 26
\\
D-53 225 BONN
\\
huebschm\@mpim-bonn.mpg.de
\endaffil
\date{October 21, 1993}
\enddate
\abstract
{Let $\Sigma$ be a closed surface, $G$ a compact Lie group, not necessarily
connected, with Lie algebra $g$, endowed with an adjoint action invariant
scalar product, let $\xi \colon P \to \Sigma$ be a principal $G$-bundle,
and pick a Riemannian metric and orientation on $\Sigma$, so that the
corresponding  Yang-Mills equations
$$d_A*K_A = 0$$
are defined, where $K_A$ refers to the curvature of a connection $A$. For
every central Yang-Mills connection $A$, the data induce a structure of
unitary representation of the stabilizer $Z_A$ on the first cohomology
group $\roman H^1_A(\Sigma,\roman{ad}(\xi))$ with coefficients in the
adjoint bundle $\roman{ad}(\xi)$, with reference to $A$, with momentum
mapping $\Theta_A$ from $\roman H^1_A(\Sigma,\roman{ad}(\xi))$ to the dual
$z^*_A$ of the Lie algebra $z_A$ of $Z_A$. We show that, for every central
Yang-Mills connection $A$,  a suitable Kuranishi map identifies a
neighborhood of zero in the  Marsden-Weinstein reduced space $\roman H_A$
for $\Theta_A$ with a neighborhood of the point $[A]$ in the moduli space
of central Yang-Mills connections on $\xi$.}
\endabstract

\endkeywords
\endsubjclass
\thanks{{\dag} The author carried out this work in the framework
of the VBAC research group of Europroj.}
\endthanks
\endtopmatter
\document
\rightheadtext{Singularities of Yang-Mills connections}
\leftheadtext{Johannes Huebschmann}

\beginsection Introduction

Let
$\Sigma$ be a closed surface,
$G$ a compact Lie group,
not necessarily connected,
with Lie algebra $g$,
and $\xi \colon P \to \Sigma$
a principal $G$-bundle.
Further, pick
a Riemannian metric on $\Sigma$ and
an {\it orthogonal structure\/}
on $g$,
that is, an adjoint action invariant positive definite
inner product.
These data then determine a Yang-Mills
theory
on the space $\Cal A(\xi)$ of connections
studied
for connected $G$
extensively by {\smc Atiyah-Bott}
in \cite\atibottw\ to which we refer for
background and notation.
In particular,
the corresponding
Yang-Mills equations look like
$$
d_A*K_A = 0
\tag0.1
$$
where $K_A$ refers to the curvature of
a connection $A$, and the solutions $A$ of (0.1) are referred to as
{\it Yang-Mills connections\/}.
Let
$\Cal N(\xi)$
be the space of
{\it  central\/} Yang-Mills connections on $\xi$,
$\Cal G(\xi)$ the group of gauge transformations, and
$N(\xi) = \Cal N(\xi)\big/\Cal G(\xi)$
the corresponding moduli space;
see Section 1 below for a precise definition of
a central connection.
In the present paper
we describe
the singularities of $N(\xi)$ explicitly.
It turns out that in a suitable sense
they are as simple as possible.
For
reasons that will become clear below
we do {\it not\/} assume
$G$ connected;
we have extended the requisite
results of
{\smc Atiyah-Bott}~\cite\atibottw\
to bundles with non-connected structure group
in~\cite\topology.
Thoughout we shall assume  that
solutions of (0.1) exist, so that the space $N(\xi)$ is non-empty.
For connected structure group this will always be so, cf. \cite\atibottw.
\smallskip
Recall that the {\it adjoint bundle\/} $\roman{ad}(\xi)$
is associated with $\xi$ via the adjoint representation
of $G$ on its Lie algebra $g$;
its sections constitute the
Lie algebra
$\bold g(\xi)$
of infinitesimal gauge transformations of $\xi$.
For a central connection $A$,
not necessarily a Yang-Mills one,
the operator
$
d_A
\colon
\Omega^*(\Sigma,\roman{ad}(\xi))
@>>>
\Omega^*(\Sigma,\roman{ad}(\xi))
$
of covariant derivative
is a differential,
that is, satisfies $d_A d_A = 0$,
even though in general
$A$ is {\it not\/}
flat;
see  (1.2) below for details.
Hence  the cohomology groups
$\roman H_A^*(\Sigma,\roman{ad}(\xi))$ are defined.
Moreover the given orthogonal structure
on $g$
induces  a symplectic structure
$
\sigma_A
$
on the
(finite dimensional)
vector space
$\roman H_A^1(\Sigma,\roman{ad}(\xi))$,
and the Lie bracket on $g$ induces a
graded Lie bracket
$
[\cdot,\cdot]_A
$
on
$\roman H_A^*(\Sigma,\roman{ad}(\xi))$
which, for degree reasons, is {\it symmetric\/}
on $\roman H_A^1(\Sigma,\roman{ad}(\xi))$.
Let $Z_A \subseteq \Cal G(\xi)$ be the stabilizer of $A$.
It is  a compact Lie group
which acts canonically on
$\roman H_A^*(\Sigma,\roman{ad}(\xi))$,
preserving
$\sigma_A$ and $[\cdot,\cdot]_A$,
and its Lie algebra $z_A$ equals
$\roman H_A^0(\Sigma,\roman{ad}(\xi))$.
The orthogonal structure
on $g$
induces a canonical isomorphism
between
$\roman H_A^2(\Sigma,\roman{ad}(\xi))$
and the dual $z_A^*$
of $z_A$
preserving the $Z_A$-actions.
Furthermore,
cf. (1.2.5) below,
the
assignment to
$\eta
\in \roman H_A^1(\Sigma,\roman{ad}(\xi))$
of
$\Theta_A(\eta) =
\frac 12 [\eta,\eta]_A$
yields a momentum mapping
$\Theta_A$
for the
$Z_A$-action
on the symplectic vector space
$\roman H_A^1(\Sigma,\roman{ad}(\xi))$,
see Section 1 below.
{\smc Marsden-Weinstein} reduction \cite\marswein\
yields the space
$\roman H_A= \Theta_A^{-1}(0)\big/ Z_A$,
and we have the following, cf. (2.32) below for a
more precise statement.

\proclaim{Theorem A}
For every central Yang-Mills connection $A$,
a suitable Kuranishi map
identifies
a neighborhood of the class $[A]$
in
$N(\xi)$
with a neighborhood of the class of zero
in
$\roman H_A$.
\endproclaim

We shall say that
a point $[A]$ of $N(\xi)$
is {\it non-singular\/}
provided
$Z_A$ acts trivially on
$\roman H^1_A(\Sigma,\roman{ad}(\xi))$.
The Theorem entails that
$N(\xi)$ is {\it smooth\/} near a
non-singular
point $[A]$
and
we have the following immediate consequence.

\proclaim{Corollary B}
The non-singular part of
$N(\xi)$
inherits from $\sigma$
a structure of
symplectic
manifold.
\endproclaim

In fact, for a central Yang-Mills connection $A$
representing a non-singular point of
$N(\xi)$,
Theorem $A$ furnishes {\it Darboux\/} coordinates
on $N(\xi)$
near the class of $A$;
in a sense,
Theorem $A$
or rather (2.32) below
yields \lq\lq Darboux coordinates\rq\rq\
near an arbitrary
point of $N(\xi)$ where
Darboux coordinates
now means the whole structure
of momentum mapping
$\Theta_A$
for the $Z_A$-action on
the symplectic vector space
$\roman H^1_A(\Sigma,\roman{ad}(\xi))$.
Theorem A makes precise the above remark that the
singularities of
$N(\xi)$
are \lq\lq as simple as possible\rq\rq.
In fact, Theorem A reduces the study of the singularities
of $N(\xi)$
to the standard example (2.4) on p.~52
of {\smc Arms-Gotay-Jennings}~\cite\armgojen.
Combining it with results
of {\smc Sjamaar-Lerman}~\cite\sjamlerm\
we shall show in \cite\singulat\
that
the decomposition
of
$N(\xi)$
into
connected components of
orbit types of central Yang-Mills connections
is a stratification
in the strong sense
in such a way that
each stratum,
being a smooth  manifold,
inherits
a finite volume symplectic structure
from the given data.
This will in general refine
the {\smc Atiyah-Bott}-decomposition
of the moduli space of {\it all\/}
Yang-Mills connections.
In particular,
for $G=U(n)$, the unitary group,
in the \lq\lq coprime
case\rq\rq, cf. {\smc Atiyah-Bott}~\cite\atibottw,
the component of the absolute minimum
of the Yang-Mills functional
has no singularities,
that is, $N(\xi)$
is smooth,
and
our
stratification
then consists of a single piece.
\smallskip
The statement of Corollary B
was obtained by
{\smc Atiyah-Bott\/}
by means of symplectic reduction
involving infinite
dimensional spaces, see p.~587 of~\cite\atibottw.
Our Corollary B avoids
this infinite dimensional symplectic reduction,
at the cost of exploiting the implicit function theorem
for Banach manifolds.
Another proof of the closedness of the
symplectic structure at the non-singular points
relying on finite dimensional techniques
has recently
been given by {\smc Weinstein}~\cite\weinstwe.
Corollary B
also paves the way
to handle arbitrary
central Yang-Mills connections
$A$, not just those yielding non-singular points
of $N(\xi)$.
The details are given
in the follow up paper
{}~\cite\direct\
where we  construct
a
{\it stratified symplectic structure\/}
in the sense
of
{\smc Sjamaar}~\cite\sjamaone\
and
{\smc Sjamaar-Lerman}~\cite\sjamlerm;
this is a Poisson structure
defined at {\it every\/}
point of
$N(\xi)$;
on each stratum, it restricts to the
corresponding symplectic Poisson structure.
As a stratified symplectic space,
for every central
Yang-Mills connection A,
the space
$\roman H_A$
will then be a local model for
$N(\xi)$
near the point represented by $A$.
Thus additional geometric information
not spelled out here is lurking behind our Theorem A,
cf. e.~g. our paper \cite\locpois.
\smallskip
In another follow up paper~\cite\singulth\
we identify the strata of $N(\xi)$
with reductions to suitable
subbundles.
Thereby we {\it cannot\/} avoid running into principal bundles
with {\it non-connected\/} structure groups,
even when the structure group of the bundle
$\xi$ we started with is connected.
This is the reason why the present theory has been set up for
general compact not necessarily connected structure groups.
\smallskip
I am indebted to
S. Axelrod,
E. Lerman,
R. Sjamaar,
P. Slodowy,
A. Weinstein,
J. Weitsman,
and T. Wurzbacher
for discussions,
and to J.~P.~Brasselet
for an illuminating lecture
about the idea of a
\lq\lq cone on a link\rq\rq\ at an early stage
of the project.

\beginsection 1. Preliminaries

{\bf 1.1. The space of connections as a K\"ahler manifold}
\smallskip\noindent
Write
$\Omega^*=\Omega^*(\Sigma,\roman{ad}(\xi))$.
The {\it data\/} we shall use
throughout
are the chosen orthogonal structure on $g$,
the Riemannian metric on $\Sigma$,
and a fixed orientation on $\Sigma$,
with unique length one volume form
$\roman{vol}_{\Sigma}$ in this orientation.
We recall from \cite\atibottw\  that the data induce
$$
\alignat1
[\cdot,\cdot]
&\colon
\Omega^*
\otimes
\Omega^*
\longrightarrow
\Omega^*,
\quad\text{a graded Lie bracket;}
\tag1.1.1
\\
\wedge &\colon
\Omega^* \otimes \Omega^*
\longrightarrow \Omega^{*}(\Sigma,\bold R),
\quad
\text{a graded commutative pairing;}
\tag1.1.2
\\
(\cdot,\cdot)&\colon
\Omega^* \otimes \Omega^{2-*}
\longrightarrow \bold R,
\quad
\text{a weakly non-degenerate
bilinear pairing,}
\tag1.1.3
\\
\intertext{given by $(\zeta, \lambda) =\int_\Sigma \zeta \wedge\lambda$;}
\cdot\,&\colon
\Omega^*
\otimes
\Omega^*
\longrightarrow \bold R,
\quad\text{a {\it weak\/} inner product;}
\tag1.1.4
\\
* &\colon \Omega^*
\longrightarrow \Omega^{2-*},
\quad\text{a duality operator.}
\tag1.1.5
\endalignat
$$
In degree one, the pairing (1.1.3) and duality operator (1.1.5)
amount to
a {\it weakly\/} symplectic structure
$\sigma= (\cdot,\cdot)$
and a complex structure
$*$ on $\Omega^1$, respectively.
\smallskip
The space
$\Cal A(\xi)$
of connections on $\xi$ is affine,
having
$\Omega^1(\Sigma,\roman{ad}(\xi))$
as its group of translations,
and hence the three pieces of structure
$\sigma,\,*,\, \cdot$
extend to
the space
of connections $\Cal A(\xi)$;
moreover they  fit together so that
$\zeta \cdot \lambda =(\zeta,*\lambda)$
whence, in particular,
they
turn  $\Cal A(\xi)$
into a K\"ahler manifold in the appropriate sense;
the K\"ahler structure is manifestly preserved by the
action of
the group
$\Cal G(\xi)$ of gauge transformations on $\Cal A(\xi)$.

\smallskip
We recall
from p.~546 of \cite\atibottw\ that
any three elements
$u,v,w$ in
$\Omega^*(\Sigma,\roman{ad}(\xi))$ satisfy
$$
[u,v]\wedge w
=
u\wedge [v,w].
\tag1.1.6
$$
This implies that, for
$|u|+|v|+|w|=2$,
$$
([u,v],w)
=
(u,[v,w]).
\tag1.1.7
$$
Next,
for every connection $A$,
when $p+q=1$,
$\phi\in \Omega^p(\Sigma,\roman{ad}(\xi))$
and
$\psi\in \Omega^q(\Sigma,\roman{ad}(\xi))$
satisfy
$$
\left(\phi,d_A\psi\right)
=(-1)^{|\phi|}\left(d_A\phi,\psi\right).
\tag 1.1.8
$$
Finally,
given $\alpha$ and $\beta$
in $\Omega^1(\Sigma,\roman{ad}(\xi))$,
we have the identity
$$
\alpha \wedge\beta =
*\alpha\wedge *\beta .
\tag1.1.9
$$
\medskip\noindent
{\bf 1.2. Central connections}
\smallskip\noindent
We denote
the centre of $G$
by
$Z$, and we write
$z$ for its Lie algebra.
It is a sub Lie algebra of the centre
of $g$,  stable under the adjoint representation.
A
{\it central\/} connection $A$ on $\xi$
is one
whose curvature $K_A \in\Omega^2(\Sigma,\roman{ad}(\xi))$
is a 2-form with values in $z$.
Thus in particular a flat connection is central.
We recall that
the topology of $\xi$
determines an element
$X_{\xi} \in z$
so that,
given an arbitrary
central
Yang-Mills connection
$A$, its
curvature $K_A$
equals
the
image
$
K_{\xi} \in \Omega^2(\Sigma,\roman{ad}(\xi))
$
of
the constant  2-form
${
X_{\xi} \otimes \roman{vol}_{\Sigma}\in
\Omega^2(\Sigma,z)
}$
under the canonical map from
$\Omega^2(\Sigma,z)$
to
$\Omega^2(\Sigma,\roman{ad}(\xi))$.
We
established
this fact in
(1.1) of our paper~\cite\topology\
for an arbitrary compact structure group,
extending a result in
Section 6 of
{\smc Atiyah-Bott\/}~\cite\atibottw\
for
a connected structure group.
Henceforth
we denote by
$
d_A
\colon
\Omega^*
@>>>
\Omega^{*+1}
$
the operator  of covariant derivative.
\smallskip
Let now $A$ be a central connection.
It is clear that
its operator of covariant derivative $d_A$
is a differential,
that is, satisfies $d_A d_A = 0$,
since $d_Ad_A = [K_A,\cdot]$.
Hence the cohomology groups
$\roman H_A^*=\roman H_A^*(\Sigma,\roman{ad}(\xi))$ are defined.
Since the operator
of covariant derivative behaves as a derivation under both the graded
Lie bracket
(1.1.1)
and the wedge product (1.1.2),
(1.1.1)~--~(1.1.3)
induce
$$
\alignat1
[\cdot,\cdot]_A
&\colon
\roman H_A^*
\otimes
\roman H_A^*
\longrightarrow
\roman H_A^*,
\quad\text{a graded Lie bracket,}
\tag1.2.1
\\
\wedge_A &\colon
\roman H_A^* \otimes \roman H_A^*
\longrightarrow \roman H^{*}(\Sigma,\bold R),
\quad
\text{a graded commutative pairing,}
\tag1.2.2
\\
(\cdot,\cdot)_A&\colon
\roman H_A^* \otimes \roman H_A^{2-*}
\longrightarrow \bold R,
\quad
\text{a non-degenerate
bilinear pairing,}
\tag1.2.3
\endalignat
$$
which, in particular,
yields the symplectic structure
$\sigma_A =(\cdot,\cdot)_A$ on $\roman H_A^1$ mentioned already
in the Introduction.
\smallskip
The kernel of the operator
$
d_A\colon
\Omega^0(\Sigma,\roman{ad}(\xi))
@>>>
\Omega^1(\Sigma,\roman{ad}(\xi))
$
of covariant derivative
is the Lie algebra $z_A$ of $Z_A$.
Hence
$\roman H_A^0(\Sigma,\roman{ad}(\xi))$
is the Lie algebra $z_A$ of $Z_A$.
Furthermore,
(1.2.3) identifies
$\roman H_A^2(\Sigma,\roman{ad}(\xi))$
with the dual
$z^*_A$ of $z_A$, and it is clear that
the
$\Cal G(\xi)$-action on
$\Cal A(\xi)$
induces an action of
$Z_A$
on
$\roman H_A^*(\Sigma,\roman{ad}(\xi))$.
Moreover, the corresponding infinitesimal
$z_A$-action on
$\roman H_A^*(\Sigma,\roman{ad}(\xi))$
is given by
the restriction of the graded Lie bracket
(1.2.1)
to
$\roman H_A^0(\Sigma,\roman{ad}(\xi))
\otimes
\roman H_A^*(\Sigma,\roman{ad}(\xi))$,
that is, by
the assignment
to $\phi \in \roman H_A^0(\Sigma,\roman{ad}(\xi))$
of the operation
$X_{\phi}\colon \roman H^*_A\to\roman H^*_A$
given by
$$
X_{\phi}(\eta) = [\phi,\eta]_A,
\quad
\eta\in \roman H^*_A(\Sigma,\roman{ad}(\xi)).
\tag1.2.4
$$
Since
the $\Cal G(\xi)$-action on
$\Cal A(\xi)$ preserves $\sigma$, it is clear that the $Z_A$-action on
$\roman H_A^1(\Sigma,\roman{ad}(\xi))$
preserves $\sigma_A$, and
we have the following,
the proof of which we leave to the reader,
cf.~\cite\golmiltw.

\proclaim{Lemma 1.2.5}
For
an arbitrary central  connection $A$,
the assignment
to
\linebreak
$\eta
\in \roman H_A^1(\Sigma,\roman{ad}(\xi))$
of
$\Theta_A(\eta) = \frac 12 [\eta,\eta]_A$
yields
a momentum
mapping
$\Theta_A$ from
$\roman H^1_A(\Sigma,\roman{ad}(\xi))$
to
$z_A^*$
for the action
of $Z_A$
on
the symplectic  vector space
$\roman H^1_A(\Sigma,\roman{ad}(\xi))$.
\endproclaim

\medskip\noindent {\bf 1.3. Hodge decomposition} \smallskip\noindent
Let $A$ be a connection on $\xi$.
As usual we
write
$
d_A^*
\colon
\Omega^*
\to
\Omega^{*-1}
$
for the adjoint of $d_A$ with
respect to the
weak inner product
(1.1.4)
on
$\Omega^*(\Sigma,\roman{ad}(\xi))$,
cf.~p.~552 of \cite\atibottw.
Since the inner product is only weak
the existence of the adjoint relies
on a suitable version
of the Hodge decomposition theorem.
\smallskip
For $j=1,2$,
we
denote by
$B^j_A(\Sigma,\roman{ad}(\xi))$
the subspace
$d_A(\Omega^{j-1}(\Sigma,\roman{ad}(\xi)))$
of coboundaries
in $\Omega^{j}(\Sigma,\roman{ad}(\xi))$
and
by
$\Cal P_A$
the orthogonal projection
from
$\Omega^j(\Sigma,\roman{ad}(\xi))$
to
$B^j_A(\Sigma,\roman{ad}(\xi))$
and, for $j=0,1,2$,
we
denote
by
$\Cal H^j_A(\Sigma,\roman{ad}(\xi))$
the vector space of {\it harmonic\/}
$j$-forms in
$\Omega^j(\Sigma,\roman{ad}(\xi))$.
The {\it Laplace\/}
operator
$$
\Delta_A = d_A d_A^*
+
d_A^*d_A
 \colon
\Omega^*(\Sigma,\roman{ad}(\xi))
@>>>
\Omega^*(\Sigma,\roman{ad}(\xi))
\tag1.3.1
$$
is manifestly $Z_A$-equivariant.
We reproduce the following
well known facts which rely on
the properties of the corresponding
{\it Green's\/} operator.

\proclaim{Proposition 1.3.2} For   a central connection $A$ on $\xi$,
for $j=1,2$,
the restriction
$$
\Delta_A| = d_A d_A^*|
 \colon
B^j_A(\Sigma,\roman{ad}(\xi))
@>>>
B^j_A(\Sigma,\roman{ad}(\xi))
\tag1.3.3
$$
is an $Z_A$-equivariant
isomorphism of real vector spaces.
\endproclaim

Let $A$ be  a central connection  on $\xi$.
For $j=0,1,2,$
we write
$$
\alpha_A \colon
\Omega^j(\Sigma,\roman{ad}(\xi))
@>>>
\Cal H^j_A(\Sigma,\roman{ad}(\xi)),\quad
\iota_A \colon
\Cal H^j_A(\Sigma,\roman{ad}(\xi))
@>>>
\Omega^j(\Sigma,\roman{ad}(\xi))
\tag1.3.4
$$
for the orthogonal projection and canonical injection,
respectively.
They are manifestly
$Z_A$-equivariant.
For $j=1,2$, we then consider the operator
$$
h_A=d_A^* \Delta_A^{-1}\Cal P_A
\colon
\Omega^j(\Sigma,\roman{ad}(\xi))
@>>>
\Omega^{j-1}(\Sigma,\roman{ad}(\xi)).
\tag1.3.5
$$
Clearly it also looks like
$
h_A = G_1 d_A^* \Cal P_A =
d_A^* G_2\Cal P_A
$
where $G_1$ and $G_2$ refer to the corresponding {\it Green's\/}
operators.
We spell out some of its properties.
\roster
\item"(1.3.6)" {\it It is  $Z_A$-equivariant\/}.
\item"(1.3.7)"
{\it It satisfies\/}
$
d^*_Ah_A =0
$
{\it and\/}
$
*h_A = d_A*\Delta_A^{-1}\Cal P_A.
$
\item"(1.3.8)" $\roman{ker}(h_A) =\roman{ker}(d^*_A)$.
\item"(1.3.9)"
For $j=1,2$, we have
$\Cal P_A =d_A h_A
\colon
\Omega^j(\Sigma,\roman{ad}(\xi))
@>>>
B_A^j(\Sigma,\roman{ad}(\xi)).
$
\endroster

The proofs of these properties are straightforward
and left to the reader.
Finally
we spell out the following
version of the
{\it Hodge decomposition theorem\/}.

\proclaim{Lemma 1.3.10}
For   a central connection $A$ on $\xi$,
the operators  $h_A$
furnish a
chain homotopy $h_A \colon  \roman{Id} \simeq \iota_A\alpha_A$,
that is, we have
$$
d_Ah_A + h_A d_A = \roman{Id} - \iota_A\alpha_A.
\tag1.3.11
$$
\endproclaim

Let now $A$ be a central connection.
Then (1.3.10) implies that
the obvious map
$$
\kappa_A
\colon
\Cal H^j_A(\Sigma,\roman{ad}(\xi))
@>>>
\roman H^j_A(\Sigma,\roman{ad}(\xi))
\tag1.3.12
$$
is an isomorphism of vector spaces; indeed, each cohomology class
in
$\roman H^j_A(\Sigma,\roman{ad}(\xi))$
has a unique harmonic representative.
Furthermore, the duality operator
(1.1.5) passes to a
duality operator
\lq\lq $\,*$\,  \rq\rq\
from $\Cal H_A^*(\Sigma,\roman{ad}(\xi))$ to
$\Cal H_A^{2-*}(\Sigma,\roman{ad}(\xi))$,
and the weak inner product (1.1.4)
induces an inner product
on
each $\Cal H_A^q(\Sigma,\roman{ad}(\xi))$
which we denote by the same symbol
\lq\lq $\,\cdot$\,  \rq\rq;
since these spaces are finite dimensional
there is no difference here between
\lq\lq weak\rq\rq\ and \lq\lq strong\rq\rq.
It is clear that these  pieces of  structure
together with the restriction of the symplectic
structure $\sigma$
constitute a hermitian structure on
$\Cal H_A^1(\Sigma,\roman{ad}(\xi))$,
having $*$ as its complex structure.
By means of the isomorphism (1.3.12), we also have this
structure
on  $\roman H_A^1(\Sigma,\roman{ad}(\xi))$;
its symplectic
structure is just $\sigma_A$.
Furthermore, the $Z_A$-action
on $\roman H_A^1(\Sigma,\roman{ad}(\xi))$
is in fact
a unitary representation since
it is compatible with all the structure, and
the momentum mapping
$\Theta_A$ is the unique one for this representation
having the value zero at the origin.

\beginsection 2. The description of the singularities

In this Section we spell out and prove a more precise version of Theorem A
in the Introduction.
Following the arguments in \cite\armamonc\
we use the Kuranishi theory of deformations
to describe the exact structure
of the singularities of our spaces of interest.
Technically it may be worthwhile remembering that,
for any connection $A$, the operator $d_A+d^*_A$ is elliptic,
and
it will be convenient to work with
suitable
Sobolev spaces.
\smallskip
The assignment
to a connection $A$
of its curvature $K_A$
is a smooth map
$J$ from
$\Cal A(\xi)$ to $\Omega^2(\Sigma,\roman{ad}(\xi))$.
It is well known to be given
by the formula
$$
J(A+\eta) =K_{A+\eta} = K_A + d_A\eta + \frac 12 [\eta,\eta],
\quad \eta \in \Omega^1(\Sigma,\roman{ad}(\xi)),
\tag2.1
$$
and its
tangent map
$
dJ(A)
$
at $A$
amounts
to the covariant derivative
$d_A$ from
$\Omega^1(\Sigma,\roman{ad}(\xi))$ to
$\Omega^2(\Sigma,\roman{ad}(\xi))$
where the tangent space
$\roman T_A\Cal A(\xi)$
is identified with $\Omega^1(\Sigma,\roman{ad}(\xi))$ as usual.
Moreover, cf. what is said in (1.2) above,
on the subspace
$\Cal N(\xi)$
of central Yang-Mills connections,
the map $J$ has a constant value
$K_{\xi}~\in~\Cal H_A^2(\Sigma,\roman{ad}(\xi))$,
determined by the topology of $\xi$.
\smallskip
Let $A \in \Cal N(\xi)$
be a smooth central Yang-Mills connection,
fixed throughout.
Since $A$ is a solution of the Yang-Mills
equations (0.1)
the value of
its curvature
$K_A = K_\xi$
lies in the space
$\Cal H_A^2(\Sigma,\roman{ad}(\xi))$
of harmonic 2-forms.
Consider the Hodge decomposition
$$
\Omega^2(\Sigma,\roman{ad}(\xi)) =
B^2_A(\Sigma,\roman{ad}(\xi)) \oplus \Cal H_A^2(\Sigma,\roman{ad}(\xi)).
\tag2.2
$$
Clearly the  point $A$ of $\Cal N(\xi)$
is {\it regular\/} for $J$ if and only if
$d_A$ is surjective,
that is, if and only if
$\roman H_A^2(\Sigma,\roman{ad}(\xi))$ is zero.
In a neighborhood of $A$,
the pre-image $J^{-1}(K_{\xi})=J^{-1}(0)$ is then a smooth manifold,
and the tangent space to $\Cal N(\xi)$ at $A$ equals
the space
$$
\roman T_A\Cal N(\xi)
=
\roman{ker}(d_A)
=
Z^1_A(\Sigma,\roman{ad}(\xi))
\tag2.3
$$
of 1-cocycles at $A$.
\smallskip
Suppose now that
$\roman H_A^2(\Sigma,\roman{ad}(\xi))$
is non-zero.
As before, let
$\Cal P_A$
be the orthogonal projection
from
$\Omega^2(\Sigma,\roman{ad}(\xi))$ onto
$B^2_A(\Sigma,\roman{ad}(\xi))$, and let
$$
\Cal A_A =
(\Cal P_A J)^{-1}(0)
= J^{-1}\left(\Cal H_A^2(\Sigma,\roman{ad}(\xi))\right).
\tag2.4
$$
Notice that
$A+\eta \in \Cal A_A$ if and only if
$
d_A^*\left(d_A\eta+\frac 12 [\eta,\eta]\right) = 0.
$
By construction,
the space $\Cal A_A$ contains
$\Cal N(\xi) $. Further,
since the composite map
$
\Cal P_A J
$
from
$\Cal A(\xi)$ to  $B^2_A(\Sigma,\roman{ad}(\xi))$
is a submersion at $A$,
in a neighborhood of $A$, $\Cal A_A$
is a smooth manifold
in such  a way that, for
$A+\eta \in \Cal A_A$,
$$
\roman T_{A+\eta}\Cal A_A
=
\{\psi; \Cal P_A d_{A+\eta} \psi=0\}
=
\{\psi; d_{A+\eta} \psi \in \Cal H_A^2(\Sigma,\roman{ad}(\xi))\}
\subseteq \Omega^1(\Sigma,\roman{ad}(\xi)).
$$
In particular,
$
\roman T_{A}\Cal A_A
=
\roman{ker}(d_A)
=
Z^1_A(\Sigma,\roman{ad}(\xi)),
$
and
$
Z^1_{A+\eta}(\Sigma,\roman{ad}(\xi))
\subseteq
\roman T_{A+\eta}\Cal A_A.
$
Define the smooth map
$$
J^{\sharp} \colon
\Cal A_A
@>>>
\Cal H_A^2(\Sigma,\roman{ad}(\xi))
\tag2.5
$$
as the restriction of $J$ so that,
for
$A+\eta \in \Cal A_A$,
$$
J^{\sharp}(A+\eta)
=
K_{\xi}+ d_A\eta+\frac 12 [\eta,\eta]
\in \Cal H^2_A(\Sigma,\roman{ad}(\xi)).
\tag2.6
$$
Note that
$J^{\sharp}(A)=K_{\xi}$ and that,
for
$A+\eta \in \Cal A_A$,
the tangent map
$dJ^{\sharp}(A+\eta)$
from
$\roman T_{A+\eta}\Cal A_A$
to
$\Cal H^2_A(\Sigma,\roman{ad}(\xi))$
is given by the restriction of the operator $d_{A+\eta}$;
in particular,
the derivative
${dJ^{\sharp}(A)
\colon Z^1_A(\Sigma,\roman{ad}(\xi)) \to
\Cal H^2_A(\Sigma,\roman{ad}(\xi))
}$
is zero
since so is the restriction of
$d_A$ to
$Z^1_A(\Sigma,\roman{ad}(\xi))$.
It is clear that the space of central Yang-Mills connections
$\Cal N(\xi)$
is smooth near $A$
and coincides with $\Cal A_A$
near $A$
if and only  if
the map $J^{\sharp}$ is constant,  having constant value $K_{\xi}$.
Hence:

\proclaim{Lemma 2.7}
For a
1-form
$\eta\in\Omega^1(\Sigma,\roman{ad}(\xi))$
having the property that
$A+\eta \in \Cal A_A$, the following are equivalent.
\roster
\item
The connection $A+\eta \in \Cal A_A$ is a central Yang-Mills
connection;
\item
$
d_A\eta + \frac 12 [\eta,\eta]=0 \in \Cal H^2_A(\Sigma,\roman{ad}(\xi));
$
\item
$
[\eta,\eta]_A=0 \in \roman H^2_A(\Sigma,\roman{ad}(\xi)). \qed
$
\endroster
\endproclaim

At this stage we thus obtain already the following well known.

\proclaim{Theorem 2.8}
The space
$\Cal N(\xi)$
of central Yang-Mills connections
coincides with $\Cal A_A$ near $A$
and hence
is smooth near $A \in \Cal N(\xi)$,
having
the space
of 1-cocycles
$Z^1_{A+\eta}(\Sigma,\roman{ad}(\xi))$
as tangent space
for every $A+\eta \in \Cal N(\xi)$
near $A$,
if and only if
the
symmetric bilinear pairing
$
[\cdot,\cdot]_A
$
on
$\roman H^1_A(\Sigma,\roman{ad}(\xi))$
is zero. \qed
\endproclaim

Recall that after a choice $Q\in \Sigma$ of base point,
the normal subgroup
$\Cal G^Q(\xi)$
of at $Q$ {\it based\/}
gauge transformations
acts freely on $\Cal A(\xi)$,
and,
after a choice $\widehat Q\in P$ of pre-image
of $Q$ has been made,
evaluation of gauge transformations at $Q$ furnishes
a surjective homomorphism
from
$\Cal G(\xi)$ onto $G$
whose kernel equals $\Cal G^Q(\xi)$.
Consequently
this
surjection
maps
the stabilizer
of an arbitrary connection
isomorphically onto a
closed
subgroup of $G$.
Under the present circumstances,
the method of Kuranishi consists of
parametrizing
a neighborhood of $A$ in $\Cal A_A$
equivariantly with respect to its stabilizer
$Z_A$
by a neighborhood of the tangent space
of $\Cal A_A$. Here are the details:
\smallskip
Recall that, in view of the Hodge decomposition theorem,
an {\it infinitesimal\/}
slice for the $\Cal G(\xi)$-action on
$\Cal A(\xi)$ is given by the
{\it transverse gauge\/},
that is, by the affine subspace
$$
\Cal S_A
= A+ \roman{ker}\left(d_A^*\colon
\Omega^1(\Sigma,\roman{ad}(\xi))
@>>>
\Omega^0(\Sigma,\roman{ad}(\xi))\right).
$$
Standard analytic arguments
involving Sobolev spaces then show that
this infinitesimial slice
generates a local slice
but for the moment this is not important for us.
We only note that
$\roman T_A\Cal S_A = \roman {ker}(d_A^*)$ and that,
with respect to (1.1.4), the Hodge decomposition
$$
\roman T_A\Cal A(\xi) =
d_A\left(\Omega^0(\Sigma,\roman{ad}(\xi))\right)
\oplus
\roman {ker}(d_A^*)
$$
is an orthogonal decomposition.
\smallskip
By means of
the operator
$h_A$, cf. (1.3.5),
the corresponding {\it Kuranishi map\/}
$F_A$ from
$\Cal A(\xi)$ to
itself
is defined by
$$
F_A(A+\eta) = A+\eta + \frac 12 h_A[\eta,\eta],
\quad
\eta\in \Omega^1(\Sigma,\roman{ad}(\xi)).
\tag2.9
$$
We shall use
its properties spelled out below,
cf. Lemmata 9~--~12
in \cite\armamonc.
\smallskip\noindent (2.10)
{\it It is\/}
$Z_A$-{\it equivariant.\/}
\newline\noindent (2.11)
{\it It is a local diffeomorphism
of a neighborhood of $A$ to a neighborhood of $A$.\/}
\newline\noindent (2.12)
{\it It satisfies the formula\/}
$$
d_A(F_A(A+\eta) - A) =  \Cal P_A(J(A+\eta)),
\quad
\eta\in \Omega^1(\Sigma,\roman{ad}(\xi)).
$$
\newline\noindent (2.13)
{\it Hence it identifies a neighborhood of\/}
$A$ {\it in\/}
$\Cal A_A$
{\it with a neighborhood of\/}
$A$ {\it in\/} $A+ Z^1_A(\Sigma,\roman{ad}(\xi))$.
\newline\noindent (2.14)
{\it For every $\eta \in \Omega^1(\Sigma,\roman{ad}(\xi))$,
we have\/}
$
d_A^*(F_A(A+\eta) - A)
=d_A^*(\eta) .
$
{\it Consequently\/}
$F_A$ {\it maps\/} $\Cal S_A$ {\it to itself
in such a way that\/}
$F_A(A+\eta) \in \Cal S_A$
{\it implies\/}
$A+\eta \in \Cal S_A$.
\smallskip
Properties (2.13) and (2.14) above imply:
\smallskip\noindent (2.15)
{\it Near $A$, the intersection
$\Cal A_A \cap \Cal S_A$
is a smooth
finite dimensional manifold,
and the Kuranishi map $F_A$ restricts to a local diffeomorphism of\/}
$\Cal A_A \cap \Cal S_A$
{\it onto\/}
$A+  \Cal H^1_A(\Sigma,\roman{ad}(\xi))$.
\smallskip
Smoothness of
$\Cal A_A \cap \Cal S_A$
near $A$ follows also from the fact that
$\Cal A_A$ and $\Cal S_A$
intersect transversely near $A$.
With respect to
(1.1.4), we now explicitly choose
a ball $\Cal B_A$ around $A$ in
$A+  \Cal H^1_A(\Sigma,\roman{ad}(\xi))$
in such a way that (i) {\it the space\/}
$$
\Cal M_A
= F_A^{-1} \left(\Cal B_A\right)
\subseteq
\Cal A_A \cap \Cal S_A,
\tag2.16
$$
{\it is a smooth finite dimensional\/}
$Z_A$-{\it manifold\/}, and (ii) {\it the Kuranishi map\/} $F_A$
{\it restricts to
a diffeomorphism\/}
$$
f_A \colon
\Cal M_A
@>>>
\Cal B_A,
\tag2.17
$$
necessarily $Z_A$-{\it equivariant\/}.
In fact,
the space
$\Cal H^1_A(\Sigma,\roman{ad}(\xi))$
being of finite dimension,
the restriction of the inner product
(1.1.4)
yields a norm
on $\Cal H^1_A(\Sigma,\roman{ad}(\xi))$
in the usual (strong) sense.
Moreover,
the action of the group
$\Cal G(\xi)$
of gauge transformations
restricts
to an action of $Z_A$
on
$\Cal H^1_A(\Sigma,\roman{ad}(\xi))$
and hence on the affine subspace
$A+  \Cal H^1_A(\Sigma,\roman{ad}(\xi))$,
and
the weak inner product (1.1.4) is invariant under gauge
transformations whence
any ball therein
is an $Z_A$-invariant subspace; in particular,
the ball
$\Cal B_A$ in $A+  \Cal
H^1_A(\Sigma,\roman{ad}(\xi))$
is $Z_A$-invariant.
Since $F_A$ is
$Z_A$-equivariant,
$\Cal M_A$ inherits
a $Z_A$-action, and
$f_A$ is $Z_A$-equivariant.
\smallskip
We now
have the machinery in place
to  show that the
diffeomorphism (2.17)
yields a symplectic change of coordinates:
Let $\omega_A$
be the restriction
of the (weakly) symplectic structure $\sigma$ on $\Cal A(\xi)$
to the smooth submanifold
$\Cal M_A$;
it is necessarily closed.

\proclaim{Lemma 2.18}
The diffeomorphism $f_A$ is compatible with the 2-forms
$\omega_A$ and $\sigma$ in the sense that, for
$A+\eta \in \Cal M_A$,
given
$
\psi,\vartheta \in \roman T_{A+\eta} \Cal M_A,
$
we have
$$
\sigma (\psi,\vartheta)
=
\omega_A (\psi,\vartheta)
=
\sigma (f_A'(A+\eta)\psi,f_A'(A+\eta)\vartheta).
$$
Consequently $\omega_A$ is non-degenerate, that is,
$\sigma$ passes to a symplectic structure on
$\Cal M_A$.
\endproclaim

\demo{Proof}
Let $A+\eta \in \Cal M_A$.
In a neighborhood of $A+\eta$,
the map $f_A$ looks like
$$
f_A(A+\eta+\psi)
=
A+\eta +
\frac 12 h_A[\eta,\eta]
+
\left(\psi+h_A[\eta,\psi]\right)
+
\frac 12 h_A[\psi,\psi],
$$
for suitable $\psi \in \Omega^1(\Sigma,\roman{ad}(\xi))$,
cf. (2.9).
Consequently its derivative
$$
f_A'(A+\eta)\colon
\roman T_{A+\eta}\Cal M_A
@>>>
\roman T_{F_A(A+\eta)}\Cal H^1_A(\Sigma,\roman{ad}(\xi))
=
\Cal H^1_A(\Sigma,\roman{ad}(\xi))
$$
at $A+\eta \in \Cal M_A$
is given by the assignment
to
$\psi \in \roman T_{A+\eta}\Cal M_A$
of
$\psi + h_A[\eta,\psi]$.
Thus the statement of the Lemma will be a consequence of the
following.
\enddemo

\proclaim{Proposition 2.19}
Given
$
\psi,\vartheta \in \roman T_{A+\eta} \Cal S_A
= \roman{ker}\left(d_A^*\colon
\Omega^1(\Sigma,\roman{ad}(\xi))
@>>>
\Omega^0(\Sigma,\roman{ad}(\xi))\right)$,
we have
$$
\sigma (\psi,\vartheta)
=
\sigma (\psi+h_A[\eta,\psi],\vartheta +h_A[\eta,\vartheta]).
$$
\endproclaim

\demo{Proof}
Clearly it will suffice to show that
$\sigma (h_A[\eta,\psi],\vartheta)$,
$\sigma (\psi,h_A[\eta,\vartheta])$,
and
\linebreak
$\sigma (h_A[\eta,\psi],h_A[\eta,\vartheta])$
are zero.
In order to see this we recall that
the duality operator $*$ is symplectic
whence
$\sigma (u,v)=\sigma (*u,*v)$,
whatever
$u,v \in \Omega^1(\Sigma,\roman{ad}(\xi))$.
By virtue of (1.3.7),
letting
$\tau
= * \Delta_A^{-1}\Cal P_A [\eta,\psi]\in \Omega^0(\Sigma,\roman{ad}(\xi))$,
we have, cf. (1.1.8),
$$
\sigma (*h_A[\eta,\psi],*\vartheta)
=
(d_A \tau,*\vartheta)
=
( \tau,d_A*\vartheta) = 0
$$
whence
$$
\sigma (h_A[\eta,\psi],\vartheta)= \sigma (*h_A[\eta,\psi],*\vartheta)= 0.
$$
The same kind of argument shows that
$\sigma (\psi,h_A[\eta,\vartheta])$ is zero.
Finally, to see that
\linebreak
$\sigma (*h_A[\eta,\psi],*h_A[\eta,\vartheta])$
is zero,
we proceed as above and obtain
$$
\sigma (*h_A[\eta,\psi],*h_A[\eta,\vartheta])
=
(\tau ,d_A*h_A[\eta,\vartheta]).
$$
However this is zero
since
$d_A*h_A$ is zero, cf. (1.3.7). \qed
\enddemo

\proclaim{Corollary 2.20}
The local diffeomorphism
$f_A$ from
$\Cal M_A$ to $A + \Cal H^1_A(\Sigma,\roman{ad}(\xi))$
is a symplectic change of coordinates,
that is, it yields Darboux coordinates on
$\Cal M_A$. \qed
\endproclaim

Our next aim is to
examine the restriction
$
J_A
$
of $J^{\sharp}$
to $\Cal M_A$,
cf. {\rm (2.5)}.
We recall
that
(1.2.3)
identifies
$\roman H_A^2(\Sigma,\roman{ad}(\xi))$
with the dual
$z^*_A$ of $z_A$
and, cf. \cite\atibottw, that
$J$, cf. (2.1),
is a  momentum mapping
for $\sigma$ and the $\Cal G(\xi)$-action
on $\Cal A(\xi)$,
the space $\Omega^2(\Sigma,\roman {ad}(\xi))$
being identified with the dual of
the Lie algebra $\bold g(\xi) =\Omega^0(\Sigma,\roman{ad}(\xi))$
of infinitesimal gauge transformations via (1.1.3).
This implies at once the following.

\proclaim{Lemma 2.21}
The composition of $J_A$
with
the canonical isomorphism
$\kappa_A$
from
$\Cal H^2_A(\Sigma,\roman{ad}(\xi))$
to
$\roman H^2_A(\Sigma,\roman{ad}(\xi))$,
cf.~{\rm (1.3.12)},
is a momentum mapping
for the action of
$Z_A$
on the symplectic manifold $\Cal M_A$. \qed
\endproclaim

\smallskip
For $\phi \in z_A=\roman H^0_A(\Sigma,\roman{ad}(\xi))$,
let
$X_{\phi}$ denote  the vector field on
$A+\Cal H^1_A(\Sigma,\roman{ad}(\xi))$
induced by the infinitesimal
$z_A$-action so that , cf. (1.2.4),
for $\eta \in \Cal H^1_A(\Sigma,\roman{ad}(\xi))$,
$$
X_{\phi}(A+\eta) = [\phi,\eta] \in \Cal H^1_A(\Sigma,\roman{ad}(\xi))
=
\roman T_{A+\eta}\left(A+\Cal H^1_A(\Sigma,\roman{ad}(\xi))\right).
\tag2.22
$$

\proclaim{Lemma 2.23}
The assignment
to
$A+\eta
\in A+\Cal H^1_A(\Sigma,\roman{ad}(\xi))$
of
$$
j_A(A+\eta) = \kappa_A(K_{\xi})+\frac 12 [\eta,\eta]_A
$$
yields a momentum
mapping $j_A$
for the action of $Z_A$
on
the affine symplectic  space
$A+\Cal H^1_A(\Sigma,\roman{ad}(\xi))$, in fact, the unique one
having
the value $\kappa_A(K_{\xi})$ at the point $A$.
\endproclaim

\demo{Proof} This is established in much the same way as (1.2.5),
by means of the canonical isomorphism (1.3.12),
combined with the canonical symplectomorphism
of affine symplectic manifolds
from
$A+\Cal H^1_A(\Sigma,\roman{ad}(\xi))$
to $\Cal H^1_A(\Sigma,\roman{ad}(\xi))$. \qed
\enddemo

As far as the statement of (2.23) is concerned,
the term involving $K_{\xi}$
may safely be ignored.
It has been included to have a consistent result in (2.24) below.

\proclaim{Theorem 2.24}
The symplectomorphism
$f_A$
preserves the momentum mappings in the sense that
the diagram
$$
\CD
\Cal M_A
@>J_A>>
\Cal H^2_A(\Sigma,\roman{ad}(\xi))
\\
@V{f_A}VV
@V{\kappa_A}VV
\\
A+\Cal H^1_A(\Sigma,\roman{ad}(\xi))
@>{j_A}>>
\roman H^2_A(\Sigma,\roman{ad}(\xi))
\endCD
$$
is commutative
where $\kappa_A$ refers to the canonical isomorphism
from
$\Cal H^2_A(\Sigma,\roman{ad}(\xi))$
to
$\roman H^2_A(\Sigma,\roman{ad}(\xi))$,
cf.~{\rm (1.3.12)}.
\endproclaim

In order to prove this we need some preparations.

\proclaim{Lemma 2.25}
For $A+\eta~\in~\Cal A_A~\cap~\Cal S_A$,
the elements
$[\eta, h_A[\eta,\eta]]$
and
$[h_A[\eta,\eta],h_A[\eta,\eta]]$
are coboundaries, that is,
pass to zero in
$\roman H^2_A(\Sigma,\roman{ad}(\xi))$.
\endproclaim

To see this we proceed as follows:
Since
$A+\eta \in \Cal A_A$
the element
$d_A\eta + \frac 12 [\eta,\eta] $
lies in $\Cal H^2_A(\Sigma,\roman{ad}(\xi))$,
and hence
$
h_Ad_A\eta + \frac 12 h_A [\eta,\eta]= 0,
$
since
$\Cal H^2_A(\Sigma,\roman{ad}(\xi)) = \roman{ker}(h_A)$,
cf.~(1.3.8).
Hence it will suffice to show that
$[\eta, h_Ad_A\eta]$
and
$[h_Ad_A\eta,h_Ad_A\eta]$
pass to  zero in
$\roman H^2_A(\Sigma,\roman{ad}(\xi))$.

\proclaim{Lemma 2.26}
For $\phi \in \roman H_A^0(\Sigma,\roman{ad}(\xi))$
and
$\eta \in \Omega^1(\Sigma,\roman{ad}(\xi))$,
$[\phi,*\eta]$ equals  $*[\phi,\eta]$.
\endproclaim

\demo{Proof}
The action of the group $\Cal G(\xi)$
of gauge transformations on $\Cal A(\xi)$
preserves $*$ in the sense that,
given a gauge transformation
$\gamma$,
we have
$$
(\roman T_A\gamma)\circ * = * \circ (\roman T_A\gamma)
\colon
\roman T_A\Cal A(\xi)
=
\Omega^1(\Sigma,\roman{ad}(\xi))
@>>>
\Omega^1(\Sigma,\roman{ad}(\xi))
=\roman T_{\gamma A}\Cal A(\xi).
$$
Consequently
the action of $Z_A$
on
$\Omega^1(\Sigma,\roman{ad}(\xi))
=\roman T_{A}\Cal A(\xi)$
preserves $*$.
Given
$\psi \in \Omega^0(\Sigma,\roman{ad}(\xi)) = \bold g(\xi)$,
let
$X_{\psi}\colon
\Cal A(\xi)
@>>>
\Omega^1(\Sigma,\roman{ad}(\xi))
$
be the vector field on
$\Cal A(\xi)$ coming from the
infinitesimal
$\bold g(\xi)$-action on $\Cal A(\xi)$;
for a connection $\widetilde A$,
it is given by
$
X_{\psi}(\widetilde A) = -d_{\widetilde A}(\psi),
$
and its derivative
$dX_{\psi}(\widetilde A)$  looks like
$$
dX_{\psi}(\widetilde A)\colon
\Omega^1(\Sigma,\roman{ad}(\xi))
@>>>
\Omega^1(\Sigma,\roman{ad}(\xi)).
$$
Since the action of the isotropy subgroup $Z_A$
on
$\Omega^1(\Sigma,\roman{ad}(\xi))
=\roman T_ A\Cal A(\xi)$
preserves the duality operator $*$,
for
$\phi \in \roman H_A^0(\Sigma,\roman{ad}(\xi)) = z_A$,
we have
$$
(\roman T*) \circ dX_{\phi}(A)
=
dX_{\phi}(A) \circ * .
$$
See Lemma 4 in \cite\armamonc\
and Proposition 4.1.28 in \cite\abramars\  for details.
Since the range of
$X_{\phi}$ is a linear space,
$\roman T*$ amounts to $*$.
However,
${
d_{A+\eta}\psi = d_A\psi + [\eta,\psi],
}$
whence
${
-dX_{\psi}(A)(\eta) = [\psi,\eta].
}$
Consequently for
$\phi \in \roman H_A^0(\Sigma,\roman{ad}(\xi)) = z_A$,
we get
$$
[\phi,*\eta] =
-dX_{\phi}(A)(*\eta)
=
-*dX_{\phi}(A)(\eta)
=
*[\phi,\eta]
$$
as asserted. \qed
\enddemo

\proclaim{Corollary 2.27}
For $\phi \in \roman H_A^0(\Sigma,\roman{ad}(\xi))$
and $\eta,\vartheta \in \Omega^1(\Sigma,\roman{ad}(\xi))$,
as real valued 2-forms,
$$
[\eta,\vartheta]\wedge \phi=
[*\eta,*\vartheta]\wedge \phi.
\tag2.27.1
$$
Consequently,
given $\eta,\vartheta \in \Omega^1(\Sigma,\roman{ad}(\xi))$,
the 2-forms
$[\eta, \vartheta]$
and
$[*\eta,*\vartheta]$
represent the same class in
$\roman H_A^2(\Sigma,\roman{ad}(\xi))$.
\endproclaim

\demo{Proof}
Applying the
identity (1.1.9)
 with $\alpha = \eta$
and
$\beta =[\vartheta,\phi]$
and keeping in mind that, in view of (2.26),
$
*[\vartheta,\phi]
=[*\vartheta,\phi]$,
we obtain, cf. (1.1.6),
$$
[\eta,\vartheta]\wedge \phi  =
 \eta \wedge [\vartheta,\phi]
=
*\eta\wedge *[\vartheta,\phi]
=
*\eta\wedge [*\vartheta,\phi]
=
[*\eta,*\vartheta]\wedge \phi.
$$
To verify the other statement,
let $\phi \in \roman H_A^0(\Sigma,\roman{ad}(\xi))$.
With reference to (1.2.3), in view of (2.27.1), we then have
$$
(\phi,[\eta, \vartheta]_A)_A=
(\phi,[*\eta,*\vartheta]_A)_A.
$$
This implies the
assertion,
since the bilinear pairing (1.2.3) is non-degenerate. \qed
\enddemo

We can now complete the proof of
Lemma 2.25.
In view of (2.27) it will suffice to show that
$[*\eta, *h_Ad_A\eta]$
and
$[*h_Ad_A\eta,*h_Ad_A\eta]$
pass to zero in
$\roman H^2_A(\Sigma,\roman{ad}(\xi))$.
However, cf. (1.3.7),
$$
*h_Ad_A\eta
=d_A* \Delta_A^{-1} d_A\eta
=d_A\psi
$$
where
$
\psi= * \Delta_A^{-1} d_A\eta
\in
\Omega^0(\Sigma,\roman{ad}(\xi)).
$
Consequently
$$
[*\eta, *h_Ad_A\eta]
=
[*\eta, d_A\psi]
=
d_A[*\eta, \psi]
-
[d_A*\eta, \psi].
$$
However, since also
$\eta \in \Cal S_A$,
$d_A^* \eta= 0$ whence  $d_A*\eta = 0$,
that is,
$[*\eta, *h_Ad_A\eta]$ equals $d_A[*\eta, \psi]$.
Likewise,
$$
[*h_Ad_A\eta, *h_Ad_A\eta]
=
[d_A\psi, d_A\psi]
=
d_A[d_A\psi, \psi].
$$
Consequently
$[*\eta, *h_Ad_A\eta]$ and
$[*h_Ad_A\eta, *h_Ad_A\eta]$ both pass to zero
in
$\roman H^2_A(\Sigma,\roman{ad}(\xi))$.
This completes the proof of (2.25). \qed
\smallskip

\demo{Proof of {\rm (2.24)}}
Let
$A+\eta \in \Cal M_A$.
In view of (2.25),
the elements
$\left[h_A[\eta,\eta],\eta\right]$
and
$\left[h_A[\eta,\eta],h_A[\eta,\eta]\right]$ pass to zero in
$\roman H^2_A(\Sigma,\roman{ad}(\xi))$.
Consequently
we have
$$
\align
j_AF_A(A+\eta)
&=
j_A\left(A+\eta+\frac 12 h_A[\eta,\eta]\right)
\\
&=
K_{\xi}+
\frac 12 \left[\eta+\frac12 h_A[\eta,\eta],\eta+\frac 12
h_A[\eta,\eta]\right]_A
\\
&=
K_{\xi}+
\frac 12 [\eta,\eta]_A
=\kappa_A f_A(A+\eta) . \qed
\endalign
$$
\enddemo

Henceforth we write
$\Cal N_A= \Cal N(\xi) \cap \Cal M_A$.

\proclaim{Corollary 2.28}
The symplectomorphism
$f_A$ maps
$\Cal N_A$
locally {\rm 1-1}
onto the cone
$$
\Cal C_A = A+
\{\eta \in \Cal H^1_A(\Sigma,\roman{ad}(\xi));
[\eta,\eta]_A = 0 \in
\roman H^2_A(\Sigma,\roman{ad}(\xi))\}.
$$
\endproclaim

Finally, let
$$
\Phi_A
\colon
\Cal M_A @>>>
\roman H^1_A(\Sigma,\roman{ad}(\xi))
\tag2.29
$$
be the smooth map
defined by
$$\Phi_A(A+\eta) = \kappa_A(F_A(A+\eta)-A)
= \kappa_A(\eta + \frac 12 h_A[\eta,\eta]),
\quad
A+\eta \in \Cal M_A,
$$
where $\kappa_A$ refers to the isomorphism (1.3.12).
In view of (2.15), this is an injective immersion,
in fact, cf. (2.20), a local symplectomorphism. Moreover, let
$$
\vartheta_A
\colon
\Cal M_A @>>>
 \roman H^2_A(\Sigma,\roman{ad}(\xi))=z_A^*
\tag2.30
$$
be the smooth map defined by
$\vartheta_A(x) = \kappa_A J_A(x) - \kappa_A(K_{\xi})$,
for $x \in \Cal M_A$.
Since $K_{\xi}$ remains invariant under the $Z_A$-action,
Lemma 2.21 implies that
$\vartheta_A$
is a momentum mapping
for the action of
$Z_A$
on the symplectic manifold $\Cal M_A$; in fact, it is the unique one
having the value zero at the point $A$.

\proclaim{Corollary 2.31}
The local diffeomorphism
$\Phi_A$ maps
$\Cal M_A$
locally 1-1
symplectically
and $Z_A$-equivariantly
onto
$\roman H^1_A(\Sigma,\roman{ad}(\xi))$
and, furthermore, preserves the momentum mappings in the sense that
$\vartheta_A = \Theta_A \Phi_A \colon \Cal M_A @>>> z_A^*$.
\endproclaim

{}From this we deduce a more precise version
of  Theorem A in the Introduction.

\proclaim{Theorem 2.32}
The local symplectomorphism
$\Phi_A$
induces a homeomorphism
of
a neighborhood of $[A]$ in
$N(\xi)$
onto a neighborhood of zero in
the Marsden-Weinstein reduced space
$\roman H_A= \Theta_A^{-1}(0)\big/ Z_A$.
\endproclaim

\demo{Proof}
Since $Z_A$ is a compact group,
a suitable Sobolev
completion of the
space $\roman{ker}(d_A^*)$
may be endowed with a $Z_A$-invariant inner product.
Any ball with respect to this inner product will then inherit
a $Z_A$-action.
By the slice theorem, cf. e.~g.~\cite\nararama,
{}~\cite\ebinone,
the map
from
$\Cal M_A\big/Z_A$
to
$\Cal A(\xi)\big/\Cal G(\xi)$
induced by the injection of $\Cal M_A$ into $\Cal A(\xi)$
is itself injective provided
the ball $\Cal B_A$
around $A + \Cal H^1_A(\Sigma,\roman{ad}(\xi))$
coming into play in (2.16)
is chosen sufficiently small;
this map restricts to a homeomorphism of
$N_A= \Cal N_A\big/Z_A$
onto a neighborhood of $[A]$
in $N(\xi)$.
On the other hand, in view of (2.31),
the
local symplectomorphism $\Phi_A$
identifies
$N_A$
locally 1-1
with $\roman H_A$. \qed
\enddemo

\smallskip
We shall say that
a central Yang-Mills connection $A$
is  {\it non-singular\/} if
its stabilizer $Z_A$ acts trivially on
$\roman H_A^1(\Sigma,\roman{ad}(\xi))$;
the point
$[A]$ of $N(\xi)$
will then be said to be {\it non-singular\/}.

\proclaim{Theorem 2.33}
Near a non-singular point
$[A]$ the space
$N(\xi)$
is smooth.
\endproclaim

\demo{Proof}
In fact, for a non-singular central Yang-Mills connection $A$
the momentum mapping $\Theta_A$ is zero
and
$\roman H_A$ coincides with
$\roman H^1_A(\Sigma,\roman{ad}(\xi))$. \qed
\enddemo

\smallskip
It may happen that
the subspace of smooth points of $N(\xi)$
is larger than that of its non-singular ones;
for example this occurs
for $G=\roman {SU}(2)$ over a surface of genus 2,
see \cite\locpois.
However the subspace
of non-singular points
is exactly that
where the symplectic structure
is defined,
that is,
the symplectic structure
on the subspace of non-singular points
cannot be extended to other smooth points.
\vfill\eject

\centerline{\smc References}
\smallskip
\widestnumber\key{999}

\ref \no  \abramars
\by R. Abraham and J. E. Marsden
\book Foundations of Mechanics
\publ Benjamin-\linebreak
Cummings Publishing Company
%\publaddr
\yr 1978
\endref

\ref \no  \armgojen
\by J. M. Arms, M. J. Gotay, and G. Jennings
\paper  Geometric and algebraic reduction for singular
momentum mappings
\jour Advances in Mathematics
\vol 79
\yr 1990
\pages  43--103
\endref

\ref \no  \armamonc
\by J. M. Arms, J. E. Marsden, and V. Moncrief
\paper  Symmetry and bifurcation of moment mappings
\jour Comm. Math. Phys.
\vol 78
\yr 1981
\pages  455--478
\endref

\ref \no  \atibottw
\by M. Atiyah and R. Bott
\paper The Yang-Mills equations over Riemann surfaces
\jour Phil. Trans. R. Soc. London  A
\vol 308
\yr 1982
\pages  523--615
\endref

\ref \no  \ebinone
\by D. Ebin
\paper The manifold of Riemannian metrics
\jour Proceedings of symposia in Pure Mathematics
\vol 15
\yr 1970
\pages  11--40
\publ American Math. Soc.
\publaddr Providence, R. I
\endref

\ref \no  \golmiltw
\by W. M. Goldman and J. Millson
\paper Differential graded Lie algebras and singularities
of level set momentum mappings
\jour Commun. Math. Phys.
\vol 131
\yr 1990
\pages 495--515
\endref

\ref \no  \singulat
\by J. Huebschmann
\paper The singularities of Yang-Mills connections
for bundles on a surface. II. The stratification
\paperinfo Math. Z. (to appear)
%\jour \vol \yr \pages
\endref

\ref \no  \singulth
\by J. Huebschmann
\paper The singularities of Yang-Mills connections
for bundles on a surface. III. The identification of the strata
\paperinfo In preparation
\endref

\ref \no  \topology
\by J. Huebschmann
\paper
Holonomies of Yang-Mills connections
for bundles on a surface with disconnected structure group
\jour Math. Proc. Cambr. Phil. Soc.
\vol 116
\yr 1994
\pages 375--384
\endref

\ref \no  \smooth
\by J. Huebschmann
\paper
Smooth structures on
certain moduli spaces
for bundles on a surface
\paperinfo preprint 1992
\endref

\ref \no  \direct
\by J. Huebschmann
\paper
Poisson structures on certain
moduli spaces
for bundles on a surface
\paperinfo Annales de l'Institut Fourier (to appear)
\endref

\ref \no \locpois
\by J. Huebschmann
\paper Poisson geometry of
flat connections for {\rm SU(2)}-bundles on surfaces
\paperinfo Math. Z. (to appear), hep-th/9312113.
%\jour \vol \yr \pages
\endref

\ref \no \modusym
\by J. Huebschmann
\paper Symplectic and Poisson structures of certain moduli spaces
\paperinfo Preprint 1993, hep-th/9312112
%\jour \vol \yr \pages
\endref

\ref \no \huebjeff
\by J. Huebschmann and L. Jeffrey
\paper Group cohomology construction
of symplectic forms on certain moduli spaces
\jour Int. Math. Research Notices
\vol 6
\yr 1994
\pages 245--249
\endref

\ref \no \marswein
\by J. Marsden and A. Weinstein
\paper Reduction of symplectic manifolds with symmetries
\jour Rep. on Math. Phys.
\vol 5
\yr 1974
\pages 121--130
\endref

\ref \no \mittvial
\by P. K. Mitter and C. M. Viallet
\paper On the bundle of connections and the gauge orbit manifold in
Yang-Mills theory
\jour Comm. in Math. Phys.
\vol 79
\yr 1981
\pages 457--472
\endref

\ref \no \nararama
\by M. S. Narasimhan and T. R. Ramadas
\paper Geometry of $\roman{SU}(2)$-gauge fields
\jour Comm. in Math. Phys.
\vol 67
\yr 1979
\pages  121--136
\endref

\ref \no \narashed
\by M. S. Narasimhan and C. S. Seshadri
\paper Stable and unitary vector bundles on a compact Riemann surface
\jour Ann. of Math.
\vol 82
\yr 1965
\pages  540--567
\endref

\ref \no  \isingeth
\by I. M. Singer
\paper Some remarks about the Gribov ambiguity
\jour Comm. in Math. Phys.
\vol 60
\yr 1978
\pages 7--12
\endref

\ref \no  \sjamaone
\by R. Sjamaar
\book Singular orbit spaces in Riemannian and Symplectic geometry
\bookinfo Thesis, University of Utrecht, 1990
\endref

\ref \no \sjamlerm
\by R. Sjamaar and E. Lerman
\paper Stratified symplectic spaces and reduction
\jour Ann. of Math.
\vol 134
\yr 1991
\pages 375--422
\endref

\ref \no \weinstwe
\by A. Weinstein
\paper On the symplectic structure of moduli space
%\jour \vol \yr\pages
\paperinfo A. Floer memorial (to appear)
\endref
\enddocument